\documentclass{article}
\usepackage{amsmath,amsfonts,amsthm}

\theoremstyle{plain} 

\newtheorem{proposition}{Proposition}[section]

\theoremstyle{definition}
\newtheorem{definition}{Definition}[section]

\theoremstyle{remark}

\newcommand{\defref}[1]{Definition~\ref{#1}}

\numberwithin{equation}{section}

\newcommand{\lie}{\mathfrak{g}}		
\newcommand{\llie}{\mathfrak{l}(\lie)}	
\newcommand{\sllie}{\mathfrak{l}^0(\lie)}
\newcommand{\lqlie}{\mathfrak{l}_h(\lie)}
\newcommand{\Lqlie}{\mathfrak{L}_h(\lie)}
\newcommand{\lqq}{\mathfrak{l}_h}	
\newcommand{\group}{G}			
\newcommand{\sll}{\mathfrak{sl}}		
\newcommand{\fun}{{\mathcal F}(\group)}	
\newcommand{\funq}{{\mathcal F}_h(\group)}
\newcommand{\funqs}{{\mathcal F}_h}	
\newcommand{\ug}{U(\lie)}		
\newcommand{\uqg}{U_h(\lie)}		
\newcommand{\C}{\mathbb{C}}		
\newcommand{\ch}{\C[[h]]}		

\newcommand{\vf}{\mathcal V}		
\newcommand{\invvf}{\vf_L}		
\newcommand{\onef}{\Gamma}		
\newcommand{\spann}[2]{#1\langle\{#2\}\rangle}	
\newcommand{\ot}{\otimes}		
\newcommand{\hot}{\otimes}		
\newcommand{\ad}{\text{ad}}		
\newcommand{\lbv}{t}			

\begin{document}

\title{
\vspace{-25mm}
\begin{flushright}\small q-alg/9608010 \\
\small KCL-TH-96-11\\[25pt]\end{flushright}
\LARGE The Problem of \\
\LARGE Differential Calculus on Quantum Groups}

\author{Gustav W. Delius\\[8pt] 
  \small Department of Mathematics, King's College London\\[-4pt]
  \small Strand, London WC2R 2LS, Great Britain\\[-4pt]
  \small e-mail: delius@mth.kcl.ac.uk\\[-4pt]
  \small http://www.mth.kcl.ac.uk/\~{}delius/}

\date{}

\maketitle

\begin{abstract}
The bicovariant differential calculi on quantum groups of Woronowicz
have the
drawback that their dimensions do not agree with that of the corresponding
classical calculus. 
In this paper we discuss the first-order differential calculus which
arises from a simple quantum Lie algebra $\lqlie$. This calculus
has the correct dimension and is shown to be bicovariant and complete. 
But it does {\em not} satisfy the Leibniz rule.
For $sl_n$ this approach leads to a differential calculus which satisfies 
a simple generalization of the Leibniz rule.

\end{abstract}
\begin{center}
Contribution to the proceedings of the 

Colloquium on Quantum Groups and Integrable Systems

Prague, June 1996
\end{center}

\section{Introduction}

In my talk at this colloquium I have reviewed the definition of quantum Lie
algebras $\lqlie$ inside the quantized enveloping algebras $\uqg$ 
\footnote{We use Drinfeld's definition of $\uqg$ as an algebra
over $\ch$.} given
in \cite{qlie1}: $\lqlie$ is an indecomposable adjoint submodule of $\uqg$  
with the same dimension as $\lie$ on which a quantum Lie product is defined
through the adjoint action. I highlighted some of the properties of
these quantum Lie algebras most of which can be found in
\cite{qlie1}-\cite{qlie5}.
Among the most interesting is the fact that their quantum Lie product 
obeys a generalization of antisymmetry which involves the $q$-conjugation
$q\leftrightarrow 1/q$. 

Instead of reiterating the contents of my talk which have already been
published, I want to use these proceedings to discuss a question which
was raised during the colloquium.
As is well known, classically the bicovariant
differential calculus on a group manifold $\group$
is determined through its Lie algebra $\lie$, 
i.e., through the space of tangent vectors at the identity of $\group$.
The question was: does a quantum Lie algebras $\lqlie$ likewise determine
a bicovariant differential calculus on a quantum group and what does
it look like? Does it shine any light on the well known problem that
a bicovariant differential calculus on a quantum group constructed
according to the framework of Woronowicz does not have the correct 
classical dimension?

Here we show how to define the first-order
differential calculus on a quantum group starting 
from our knowledge of its quantum Lie algebra $\lqlie$.
The resulting calculus automatically has the right dimension, is complete
and bicovariant.
Unlike Woronowicz we do not impose the Leibniz rule. This is in the
spirit of Faddeev's old idea that the Leibniz rule should be
modified. In our approach the action of the differential operator $d$ on
a product is something which can be calculated and 
it is indeed found {\em{not}}
to satisfy the standard Leibniz rule. Rather it depends on the form which
the $\uqg$ coproduct takes on $\lqlie$.

For $\sll_n$ Sudbery \cite{Sud} has shown that one can choose a particular
$\lqq(\sll_n)$ on which the $\uqg$ coproduct takes the rather
simple form \eqref{sudco}. The differential calculus which we construct from
this quantum Lie algebra satisfies a very simple generalization 
\eqref{sudleibniz}
of the Leibniz rule while of course maintaining the correct dimension.
A calculus of the correct dimension with a modified Leibniz rule has 
also been proposed by Faddeev and Pyatov \cite{Fad}. The relation between
these calculi should be studied. It is possible
that they agree which would be interesting because
our construction is so very simple and straightforward.

We expect that the results of this paper on the
existence and bicovariance of the first-order differential calculus
arising from our quantum Lie algebras will seem trivial to the
experts of the field. 
We nevertheless hope that this modest 
contribution will be stimulating to the reader.

\section{Preliminaries on quantum groups}

In the theory of quantum groups one does not
actually quantize the group manifold $\group$ but rather the algebra 
$\fun$ of regular functions on the group manifold. Classically $\fun$
is the Hopf dual of $\ug$, the universal enveloping algebra of the Lie
algebra $\lie$ of $\group$. This means that the algebra of functions is
spanned by the matrix elements of all the finite-dimensional irreducible
representations of $\ug$. We define the quantized algebra of
functions $\funq$ similarly as the Hopf dual of the quantized enveloping
algebra $\uqg$. Here $\uqg$ is the algebra over $\ch$ introduced by
Drinfel'd.

\begin{definition}\label{funqdef}
Let $\group$ be a connected, simply-connected, finite-dimensional,
complex simple Lie group with Lie algebra $\lie$. The
{\em quantized algebra of functions} $\funq$ on $\group$ is the
Hopf algebra dual of the quantized enveloping algebra $\uqg$ in
the following sense.
Let $\{V^\mu\}_{\mu\in D_+}$ be the set of all finite-dimensional 
indecomposable $\uqg$-modules.
Choose bases $\{v_\mu^a\}_{a=1,\cdots,\text{dim} V^\mu}$ 
for $V^\mu$ and introduce representation functionals
$\pi^{\mu}_a{}^b:\uqg\rightarrow\ch$ by\footnote{We use the summation 
convention according to which repeated indices are summed over.}
\begin{equation}
x\,v_\mu^b=v_\mu^a\,\pi^{\mu}_a{}^b(x),~~~~
\forall x\in\uqg.
\end{equation}
These functionals span $\funq$,
\begin{equation}
\funq=\spann{\ch}{\pi^{\mu}_a{}^b}_{\mu\in D_+; a,b=1\dots \text{dim}V^\mu}.
\end{equation}
The coalgebra structure of $\funq$ is given by the algebra
structure of $\uqg$,
\begin{gather}\label{coprod}
\Delta(\pi^\mu_a{}^b)(x\hot y)\equiv
\pi^\mu_a{}^b(x\cdot y)
=(\pi^\mu_a{}^c\hot\pi^\mu_c{}^b)(x\hot y),
\\
\epsilon(\pi^\mu_a{}^b)\equiv\pi^\mu_a{}^b(I)=\delta_a{}^b.
\end{gather}
Similarly the algebra structure of
$\funq$ is determined by the coalgebra structure of $\uqg$
leading to an expression of the 
product on $\funq$ in terms of the $\uqg$ Clebsch-
Gordan coefficients,
\begin{equation}\begin{split}\label{intertwine}
(\pi^\mu_a{}^b\cdot\pi^\nu_c{}^d)(x)&\equiv
(\pi^\mu_a{}^b\hot\pi^\nu_c{}^d)(\Delta(x))\\
&=\sum_\lambda(K^{\mu\nu}_\lambda)_{ac}{}^i(K_{\mu\nu}^\lambda)_j{}^{bd}~
\pi^\lambda_i{}^j(x),
\end{split}\end{equation}
Here $K_{\mu\nu}^\lambda$ is the matrix of Clebsch-Gordan coefficients
for the decompsition of $V^\mu\otimes V^\nu$ into $V^\lambda$ and
$K^{\mu\nu}_\lambda$ is its inverse. The second equality above is
simply the defining relation for Clebsch-Gordan coefficients.
We suppress the multiplicity
label which is necessary if the embedding of $V^\lambda$ into 
$V^\mu\otimes V^\nu$ is not unique.
\end{definition}


Equation \eqref{intertwine} also shows that in those cases where all 
irreducible representations $\pi^\mu$ can be obtained from the multiple
tensor product of one fundamental representation 
$\pi^{\text{fund}}\equiv T$, the whole
algebra $\funq$ is generated by the $T_a{}^b$ alone. Relations
among these generators are obtained from the universal R-matrix 
$R$ of $\uqg$,
\begin{align}\label{RTT}
(T_a{}^b\cdot T_c{}^d)(x)&=(T_a{}^b\hot T_c{}^d)(\Delta(x))\notag\\
&=(T_a{}^b\hot T_c{}^d)(R^{-1}\Delta^{\text{op}}(x)\,R)\\
&=((R^{-1})_{ac}{}^{a'c'}T_{c'}{}^{d'}\cdot T_{a'}{}^{b'}
R_{b'd'}{}^{bd})(x)\notag,
\end{align}
where $R_{ab}{}^{cd}=(T_a{}^c\ot T_b{}^d)R$, etc.
Additional relations express the linear dependence among the $T_a{}^b$.
These are the relations of the FRT algebras of Faddeev, Reshetikhin
and Takhtajan \cite{FRT}. The FRT algebras are defined only for $\group=SL_n(\C),SO_n(\C)$ and $Sp_n(\C)$. 
It is believed that for $\group=SL_n(\C)$ and $Sp_n(\C)$ the FRT
algebras coincide with $\funq$ defined in Definition \ref{funqdef}.
This can not be the case for the groups $SO_n(\C)$ because they are not 
simply connected. Because in this paper
we are interested in constructing the differential geometry on the
group manifold by group translation from the tangent vector at the
origin, we work with the quantized function algebra on their simply connected covering groups as defined in \defref{funqdef}
rather than with the FRT algebras.

\section{Quantum Lie algebras}

We will be very brief here, see \cite{qlie1}-\cite{qlie5} for details.

A Lie algebra $\lie$ carries the adjoint representation of $\ug$. Let
$\pi^\Psi$ denote the adjoint representation of $\uqg$
(this is the unique indecomposable representation of $\uqg$ with
dimension equal to $\dim\lie$ and is a deformation of the classical
adjoint representation). 

\begin{definition}
A {\em weak quantum Lie algebra} $\lqlie$ is a submodule of $\uqg$ 
which generates $\uqg$ and on which the adjoint action gives the
adjoint representation.
The {\em quantum Lie product} is given by 
$[a,b]_h=(\ad\,a)\, b$ for all $a,b\in\lqlie$.
\end{definition}

This is a reformulation of the definition given in \cite{qlie1}
making use of the understanding gained in \cite{qlie3}. 
{\em Strong} quantum Lie algebras
$\Lqlie$ have additional properties regarding the action of the
antipode, the Cartan involution and any diagram automorphism, see
\cite{qlie1}.

By the statement that the adjoint action on $\lqlie$ gives the adjoint
representation we mean that $\lqlie$ has a
basis $\{\lbv^i\}_{i=1,\cdots\dim\lie}$ which satisfies
\begin{equation}
(\ad\, x)\,\lbv^i=\lbv^j\,\pi^\Psi_j{}^i(x),~~~~~~\forall x\in\uqg.
\end{equation}
When we say that $\lqlie$ generates $\uqg$ we mean that any element
of $\uqg$ can be expressed as a formal power series in $h$ whose
coefficients are polynomials in the $\lbv^i$.

The representation matrices of the quantum Lie algebra generators $\lbv^i$
are given by Clebsch-Gordan coefficients,
\begin{equation}\label{rep}
\pi^\mu_a{}^b(\lbv^i)=(K_{\mu\Psi}^\mu)_a{}^{bi},~~~~~
\forall \mu\in D_+,
\end{equation}
This can be derived by applying $\pi^\mu_a{}^b$ to the 
equality\footnote{Here we use the definition of the adjoint action 
$(\ad x)y=x_{(2)}\,y\,S^{-1}(x_{(1)})$.}
$x\,\lbv^i=((\ad x_{(2)})\,\lbv^i)\,x_{(1)}$ and
recognizing the result as the intertwining property 
of the Clebsch-Gordan 
coefficients. This implies in particular that the structure constants
of the quantum Lie algebra are given by the Clebsch-Gordan 
coefficients for adjoint $\ot$ adjoint into adjoint,
\begin{equation}\label{st}
[\lbv^i,\lbv^j]=(K_{\Psi\Psi}^\Psi)_k{}^{ji}\,\lbv^k.
\end{equation}
For more details see \cite{qlie3}.

Inside any $\uqg$ there exist infinitely many quantum Lie algebras
$\lqlie$, see \cite{qlie3} for their construction. However the equations
$\eqref{rep},\eqref{st}$ hold for any of them, the only differences
being in the choice of $K$. In particular,
except for $\lie=\sll_n$ with $n>2$,
all quantum Lie algebras $\lqlie$ inside 
$\uqg$ are isomorphic. This is due to the uniqueness of the embedding
of the adjoint representation into adjoint $\ot$ adjoint. In the
case of $\lie\neq\sll_n$ with $n>2$ uniqueness holds after imposing
the extra requirement of invariance under the diagram automorphism
of $U_h(\sll_n)$. Because all $\lqlie$ are isomorphic we view them
as different embeddings of an abstract quantum Lie algebra $\lie_h$
into $\uqg$.

The fact that there exist infinitely many isomorphic quantum Lie algebras 
inside $\uqg$ is not
surprising. Also classically there are infinitely many embeddings
of a Lie algebra $\lie$ into the enveloping algebra $\ug$
(multiplying all elements of one embedding by a central element of
$\uqg$ gives a new embedding).
However, among these there is a unique embedding, which
we will denote by $\sllie$, with the
property that $\sllie$ generates $\ug$. 
(This is the natural embedding and therefore commonly $\sllie$ would be 
denoted simply as $\lie$). The requirement that
$\lqlie$ generate $\uqg$ implies that $\lqlie=\sllie\mod h$. However
this requirements does not put any restrictions on $\lqlie$ at higher
orders in $h$.

Classically $\sllie$ is also special among all embeddings $\llie$
because on $\sllie$ the coproduct of $\ug$ takes the particularly 
simple form
\begin{equation}\label{special}
\Delta(a)=1\ot a+a\ot 1,~~~\forall a\in\sllie.
\end{equation}
We expect that also at the quantum level some $\lqlie$ will be 
singled out by their property under the $\uqg$ coproduct. This property
should be some quantum generalization of \eqref{special}. Unfortunately
the correct generalization has so far not been found in general. 
For $\lie=\sll_n$ Sudbery \cite{Sud} has shown that one can choose a
$\lqq^0$ which satisfies the following generalization of \eqref{special}
\begin{equation}\label{sudco}
\Delta(\lbv^i)=C\ot\lbv^i+\lbv^j\ot f_j{}^i,
\end{equation}
where $f^j_i$ are certain elements in $\uqg$ and $C$ is a central
element of $\uqg$. We will see the 
relevance of this for the Leibniz rule 
in the quantum differential calculus later.

\section{First-order differential calculus}
 
For background on differential calculus on quantum groups I recommend
the pioneering paper by Woronowicz \cite{Wor} and the review \cite{Asc}.
We take vector fields to be the starting point rather than forms.
Such an approach has been studied in \cite{vec}.

A quantum Lie algebra $\lqlie$ is interpreted as the quantum analog
of the space of tangent vectors to the group manifold at the identity.
By left-translating a tangent vector $\lbv$ one obtains a left-invariant
vector field $\lbv_L$. In order to be applied at the quantum level, the
notion of left-translation has to be formulated in terms of the functions
on the group manifold rather than the points. As explained in
\cite{Wor,Asc} this leads to
\begin{equation}\label{leftvec}
\lbv_L(a)=a_{(1)}\,\lbv(a_{(2)}),~~~~~~~~
a\in\funq,\ \lbv\in\lqlie.
\end{equation}
Here the action of $t$ on $a_{(2)}$ is given by the duality between
$\uqg$ and $\funq$, i.e., $\lbv(a_{(2)})=a_{(2)}(t)$.
In this way the quantum Lie algebra generators $\lbv^i$ lead to basis
vectors $\lbv_L^i$ for the space\footnote{
Here we use the word "space" even
when the base is only a ring, the term "free module" would be more 
precise.} 
of left-invariant vector fields,
$\invvf=\spann{\ch}{\lbv_L^i}$.

Classically an arbitrary vector field can be written as a linear 
combination of products of
left-invariant vector fields by functions. In
the quantum case we define the product $(t_L a)$ of a left-invariant
vector fields $t_L$ by an element $a$ of $\funqs$ from the right by
$(\lbv_L a)(b)=t_L(b)\,a$.
\begin{definition}
The space
$\vf$ of {\em vector fields} on the quantum group is
the right $\funqs$-module freely generated by $\{\lbv^i_L\}$, 
$\vf=\langle\{\lbv_L^i\}\rangle\funqs$.
\end{definition}
The reason why we choose to multiply vector fields by functions from the 
right is that then one-forms will be multiplied from the left, as we will
see below. Classically there is of course no distinction.

We could alternatively start with the right invariant-vector fields 
$\lbv^R$ obtained by right translating the Lie algebra elements,
giving $\lbv^R(a)=\lbv(a_{(1)})\,a_{(2)},~
a\in\funq,\ \lbv\in\lqlie$. We could use these to define the space of all
vector fields as $\vf=\langle\{\lbv_R^i\}\rangle\funqs$. Our formalism
would not be natural if this lead to a different notion of vector
fields. Luckily we have

\begin{proposition}
Left- and right-invariant vector fields are related by
\begin{equation}\label{re}
\lbv_L^i=\lbv^j_R\,\pi^\Psi_j{}^i.
\end{equation}
\end{proposition}
\begin{proof}
We need to show that
$\lbv^i_L(\pi^\mu_a{}^c)=\lbv_R^j
(\pi^\mu_a{}^c)\,\pi^\Psi_j{}^i,~\forall\pi^\mu_a{}^c$. 
Using the definitions of $t^i_L$ \eqref{leftvec} and $t^i_R$ and the 
formula \eqref{coprod} for the coproduct this becomes
$\pi^\mu_a{}^b\lbv^i(\pi^\mu_b{}^c)=
t^j(\pi^\mu_a{}^b)\,\pi^\mu_b{}^c\,\pi^\Psi_j{}^i$ 
which holds because
of \eqref{rep} and the intertwining property of the Clebsch-Gordan
coefficients.
\end{proof}
Conversely we have $\lbv_R^i=\lbv^j_L\,((\pi^\Psi)^{-1})_j{}^i$
and as a consequence $\vf=\langle\{\lbv_L^i\}\rangle\funqs=
\langle\{\lbv_R^i\}\rangle\funqs$. We note that it was not a priori obvious 
that this would work. In order to have the same dimension as
classically, we have dropped some of the
left-invariant vector fields $t^i_L$ which appear in the Woronowicz calculus. It could have 
happened that the expression of the right-invariant vector fields in terms
of the left-invariant ones would have reintroduced the ones we had dropped
and our truncation would not have been consistent.

We now introduce the space $\onef_L$ of left-invariant one-forms as the dual
space to the space $\invvf$ of left-invariant vector fields and we choose
a basis $\{\omega_i^L\}$ such that $\omega_i^L(
\lbv^j_L)=\delta_i{}^j$.
We define the action of a $\omega^L\in\onef$ on an arbitrary vector field
by setting
$\omega^L(t_L\,a)=(\omega^L(t_L))a,~\forall t_L\in\invvf, a\in\funqs$.
We define the left multiplication of one-forms with functions by
$(a\,\omega)(t)=a(\omega(t)),~\forall t\in\vf, a\in\funqs$.
\begin{definition}
The space
$\onef$ of {\em one-forms} on the quantum group is
the left $\funqs$-module freely generated by $\{\omega_i^L\}$, 
$\onef=\funqs\langle\{\omega_i^L\}\rangle$.
\end{definition}
Again we could equivalently have used right-invariant one-forms
$\omega_i^R$ satisfying $\omega_i^R(\lbv^j_R)=\delta_i{}^j$.
These are related to the left-invariant $\omega_i^L$ by
\begin{equation}
\omega^R_i=\pi^\Psi_i{}^j\omega^L_j,
\end{equation}
and thus $\onef=\funqs\langle\{\omega_i^L\}\rangle=
\funqs\langle\{\omega_i^R\}\rangle$.

As explained in \cite{Wor,Asc}, 
left- and right-translation of one-forms is described by
maps $\Delta_L:\onef\rightarrow\funqs\otimes\onef$ and
$\Delta_R:\onef\rightarrow\onef\otimes\funqs$. 
They act as
\begin{equation}
\Delta_L(a\,\omega^L)=\Delta(a)\,(1\otimes\omega^L),~~~~
\Delta_R(a\,\omega^R)=\Delta(a)\,(\omega^R\otimes 1).
\end{equation}
Because any one-form can be expressed in terms of left- or
right-invariant forms, this determines left- and right-translation
entirely.
Because of the way in which we have defined the left- and right-
translations the bicovariance condition
$(1\otimes\Delta_R)\Delta_L=(\Delta_L\otimes 1)\Delta_R$ is
guaranteed. Of course it can also be checked by explicit calculation.

\begin{definition}\label{ddef}
The {\em exterior differential} $d:\funqs\to\onef$ is defined by
\begin{equation}
da(v)=v(a)~~~~\forall a\in\funqs,\,v\in\vf.
\end{equation}
\end{definition}
In terms of the bases this can be expressed as
\begin{equation}\label{dif}
da=\lbv^i_L(a)\,\omega_i^L=\lbv^i_R(a)\,\omega_i^R.
\end{equation}

\begin{proposition}\label{p}
The first-order differential calculus $(\onef,d)$ satisfies
\begin{enumerate}
\renewcommand{\labelenumi}{(\roman{enumi})}
\item
Completeness:  
any $\rho\in\onef$ can be written as
$\rho=\sum_k\,a_k\,d\,b_k$ for some $a_k,b_k\in\funqs$.
\vspace{-1mm}
\item
Bicovariance: for all $a,b\in\funqs$,
\begin{equation}
\Delta_L(a\,db)=\Delta(a)\,(1\otimes d)\,\Delta(b),~~~~
\Delta_R(a\,db)=\Delta(a)\,(d\otimes 1)\,\Delta(b).
\end{equation}
\end{enumerate}
\end{proposition}
\begin{proof}
i) Here we can follow \cite{Asc}.
We will show that any $\omega_L^i$ can be written in this form.
Define the matrix $\gamma^{ij}=\pi^{\Psi}_a{}^b(\lbv^i)\,
\pi^{\Psi}_b{}^a(\lbv^j)$. Classically $\gamma$ goes over into the
Killing metric and thus $\gamma$ is invertible also in the
quantum case. It can now be checked that\footnote{
Remember that $t^i(\pi^{\Psi b}_a)\equiv \pi^{\Psi b}_a(t^i)$.}
\begin{equation}
\omega^L_i=(\gamma^{-1})_{ij}\pi^\Psi_a{}^b(t^j)((\pi^\Psi)^{-1})_b{}^c
\,d\,\pi^\Psi_c{}^a.
\end{equation}
ii)by a simple calculation: $\Delta_L(a\,db)=\Delta(a)
\Delta_L(t^i_L(b)\omega_i^L)=
\Delta(a)\Delta(b_{(1)})t^i(b_{(2)})(1\otimes\omega_i^L)=
\Delta(a)(b_{(1)}\otimes b_{(2)} t^i(b_{(3)})\omega^L_i)=
\Delta(a)(b_{(1)}\otimes db_{(2)})=\Delta(a)(1\otimes d)
\Delta(b)$.
\end{proof}

Woronowicz, in his approach to bicovariant quantum differential calculus,
postulates completeness and bicovariance as axioms rather than 
deriving them. He postulates one more axiom: the Leibniz rule
\begin{equation}\label{jacobi}
d(ab)=a(db)+(da)b,~~~~~\forall\,a,b\in\funqs
\end{equation}
In the second term this equation involves multiplication of
a one-form by a function from the right. 
It thus makes sense only if $\onef$ is
a $\funqs$-bimodule. The Leibniz rule is the only reason why one
is interested in having a bimodule structure on $\onef$. It is
introduced by defining the
right multiplication of one-forms by functions through
\begin{equation}\label{rightmult}
\omega_i^L\,a=a_{(1)}\,f_i{}^j(a_{(2)})\,
\omega^L_j,~~~~~\forall\,a\in\funqs,
\end{equation}
where $f^i{}_j$ are certain elements of $\uqg$.

In our approach the Leibniz rule is not an axiom but rather the
action of $d$ on a product of functions can be computed in terms
of the known coproduct $\Delta(t^i)=t^i_{(1)}\otimes t^i_{(2)}$.
\begin{equation}
d(ab)=a_{(1)}b_{(1)}t^i(a_{(2)}b_{(2)})\omega^L_i=
a_{(1)}b_{(1)}t^i_{(1)}(a_{(2)})t^i_{(2)}(b_{(2)})\omega^L_i.
\end{equation}
Even though all quantum Lie algebras $\lqlie$ inside the same $\uqg$
are isomorphic as algebras, the coproducts $\Delta(t^i)$ differ and
can be very complicated in general. 
It is important to define the differential calculus in terms of 
particular quantum Lie algebras $\lqq^0(\lie)$ whose coproduct
is simple, otherwise one will obtain an unmanageable Leibniz rule.
For $\sll_n$ a good $\lqq^0(\sll_n)$ has been found by Sudbery
\cite{Sud} which has the coproduct given in \eqref{sudco}. It leads
to a differential calculus with generalized Leibniz rule
\begin{equation}\label{sudleibniz}
d(ab)=c(a)(db)+(da)b,~~~~~\forall\,a,b\in\funqs,
\end{equation}
where $c(a)\equiv a_{(1)}C(a_{(2)})$ and right multiplication is
defined as in \eqref{rightmult}.

\section{Discussion}

We have explained that by starting from the quantum Lie algebras $\lqlie$
as studied in \cite{qlie1}-\cite{qlie5} and \cite{Sud} one automatically 
obtains complete, bicovariant quantum differential
calculi $(\onef,d)$ of the correct dimension which do however not satisfy
the standard Leibniz rule. The straightforward construction can be
summarized as
follows: 1) interpret the quantum Lie algebra elements as tangent vectors
at the identity, 2) left- or right-translate them to obtain left- or
right-invariant vector fields, 3) take their duals to obtain left- or
right invariant one-forms, 4) multiply them by functions to obtain
the space $\onef$ of all one-forms, 5) define the exterior differential
as in \defref{ddef}.

This paper has been kept very brief. Many properties of the differential
calculus have been left unstudied, in particular the exterior calculus
and the Cartan calculus. For the later we
would like to draw the reader's attention to the elegant work of 
Schupp and coworkers, see e.g. \cite{Schupp}.

The big challenge which remains is to identify among the infinitely many
quantum Lie algebras $\lqlie$ those with the simplest coproduct which
will lead to the differential calculus with the most natural generalization
of Leibniz rule. This has so far been accomplished only for $\sll_n$ by
Sudbery \cite{Sud}.

For an extensive list of references on quantum Lie algebras and
on differential calculus on quantum groups visit the World Wide
Web site on quantum Lie algebras at 
http://www.mth.kcl.ac.uk/\~{}delius/q-lie.html


\begin{thebibliography}{9}

\bibitem{Asc}P. Aschieri, L. Castellani, {\it An introduction to
noncommutative differential geometry on quantum groups}, Int. J. Mod. Phys.
{\bf A8} (1993) 1667, hep-th/9207084.

\bibitem{vec}P. Aschieri, P. Schupp, {\it Vector Fields on Quantum
Groups}, Int. J. Mod. Phys. A11 (1996) 1077-1100. q-alg/9505023.

\bibitem{qlie1} G.W. Delius, A. H\"uffmann, {\it On Quantum Lie Algebras
 and Quantum Root Systems}, q-alg/9506017, 
 J. Phys. A {\bf 29} (1996) 1703.

\bibitem {qlie2} G.W. Delius, M.D. Gould, A. H\"{u}ffmann, Y.-Z. Zhang,
 {\it Quantum Lie algebras associated to $U_q(gl_n)$ and $U_q(sl_n)$},
 J. Phys. A., q-alg/9508013.

\bibitem{qlie3} G.W. Delius, M.D. Gould, {\it Quantum Lie algebras, their
existence, uniqueness and q-antisymmetry}, q-alg/9605025.

\bibitem{qlie4} G.W. Delius, {\it Introduction to Quantum Lie Algebras},
q-alg/9605026.

\bibitem{qlie5} G.W. Delius, C. Gardner, {\it The structure of the quantum
Lie algebras of type $B_n, C_n$ and $D_n$}, KCL-TH-96-12, in preparation.

\bibitem{Fad} L.D. Faddeev, P.N. Pyatov, {\it
The differential calculus on quantum linear groups}, hep-th/9402070.

\bibitem{FRT} L.D. Faddeev, N.Yu. Reshetikhin, L.A. Takhtajan,
Algebra and Analysis {\bf 1} (1987) 178.

\bibitem{Schupp} P. Schupp, {\it Cartan Calculus: Differential Geometry
for Quantum Groups}, Proceedings of the International 
School of Physics "Enrico Fermi" Course CXXVII (1994), Editors 
L. Castellani and J. Wess, IOS Press, 1996. hep-th/9408170.

\bibitem{Sud} A. Sudbery, {\it Quantum Lie algebras of type $A_n$},
q-alg/9510004.

\bibitem{Wor} S.L. Woronowicz, {\it Differential Calculus on Compact
Matrix Pseudogroups (Quantum Groups)}, Commun. Math. Phys. {\bf 122}
(1989) 125.

\end{thebibliography}
\end{document}